\documentclass[11pt,reqno]{article}
\usepackage{graphicx}
\usepackage{amsmath}
\usepackage{amsthm}
\usepackage{latexsym,bm}
\usepackage{amssymb}
\usepackage{indentfirst}
\newtheorem{thm}{Theorem}[section]
\newtheorem{cor}[thm]{Corollary}

\numberwithin{equation}{section}

\begin{document}

\title{{\bf  Heinz type estimates for graphs in Euclidean space
}}

\author{Francisco Fontenele\thanks{Partially supported by CNPq (Brazil)}\\ \\\small{\it To my wife
Andrea}}

\date{}
\maketitle

\begin{quote}
\small {\bf Abstract}. Let $M^n$ be an entire graph in the
Euclidean $(n+1)$-space $\mathbb R^{n+1}$. Denote by $H$, $R$ and
$|A|$, respectively, the mean curvature, the scalar curvature and
the length of the second fundamental form of $M^n$. We prove that
if the mean curvature $H$ of $M^n$ is bounded then $\inf_M|R|=0$,
improving results of Elbert and Hasanis-Vlachos. We also prove
that if the Ricci curvature of $M^n$ is negative then
$\inf_M|A|=0$. The latter improves a result of Chern as well as
gives a partial answer to a question raised by Smith-Xavier. Our
technique is to estimate $\inf|H|,\;\inf|R|$ and $\inf|A|$ for
graphs in $\mathbb R^{n+1}$ of $C^2$-real valued functions defined
on closed balls in $\mathbb R^n$.



\end{quote}

\section{Introduction}

Let $B_r\subset\mathbb R^n$ be an open ball of radius $r$, and
$f:B_r\to\mathbb R$ a $C^2$-function. For $n=2$, Heinz \cite{H}
obtained the following estimates for the mean curvature $H$ and
the Gaussian curvature $K$ of the graph of $f$:

\begin{eqnarray}\label{HH}
\inf |H|\leq\frac{1}{r},
\end{eqnarray}

\begin{eqnarray}\label{HK}
\inf |K|\leq\frac{3\,e^2}{r^2},
\end{eqnarray}
where $e$ is the basis for the natural logarithm. Chern \cite{C}
and Flanders \cite{Fl}, independently, obtained inequality
(\ref{HH}) for any $n\geq 2$. As an immediate consequence one has
that an entire graph in $\mathbb R^{n+1}$, i.e., the graph of a
$C^2$-function from $\mathbb R^n$ to $\mathbb R$, cannot have mean
curvature bounded away from zero. This result of Chern and
Flanders implies that an entire graph in $\mathbb R^{n+1}$ with
constant scalar curvature $R\geq 0$ satisfies $R=0$ (see Section
2, equality (\ref{HRA})).

\vskip5pt

Salavessa \cite{Sa} extended inequality (\ref{HH}) to graphs of
smooth real valued functions defined on oriented compact domains
of Riemannian manifolds, and Barbosa-Bessa-Montenegro \cite{BBM}
extended it to transversally oriented codimension one
$C^2$-foliations of Riemannian manifolds.

\vskip5pt

In the theorems below, by a graph over a closed ball $\overline
B_r\subset\mathbb R^n$ we mean the graph in $\mathbb R^{n+1}$ of a
$C^2$-real valued function defined in $\overline B_r$. In our
first result, we improve the estimates given by Heinz, Chern and
Flanders for $\inf|H|$, by showing that the estimates can be made
strict if we consider graphs over closed balls instead of graphs
over open balls.

\begin{thm}\label{FH} If $M^n\subset\mathbb R^{n+1}$ is a graph over a closed ball $\overline B_r\subset\mathbb R^n$,then
\begin{eqnarray}\label{FH1}
\inf_M|H|<\frac{1}{r}.
\end{eqnarray}
\end{thm}

\vskip10pt



The estimate (\ref{HK}) implies that an entire graph in $\mathbb
R^3$ cannot have negative Gaussian curvature bounded away from
zero, a result extended later to complete surfaces in $\mathbb
R^3$ by Efimov \cite{E} in a remarkable work (see the discussion
after Corollary \ref{E-HV}). In the next result, we obtain a
version in higher dimensions of the inequality (\ref{HK}).

\vskip10pt

\begin{thm}\label{FR} Let $M^n\subset\mathbb R^{n+1}$ be a graph over a closed ball $\overline B_r\subset\mathbb R^n$, and denote by
$R$ the scalar curvature of $M$. Then
\begin{eqnarray}\label{FR1}
\inf_M|R|\leq\Big(\sup_M|H|+\frac{1}{r}\Big)\frac{2}{r}.
\end{eqnarray}
Moreover, if $M$ has a point where the second fundamental form is
semi-definite, then
\begin{eqnarray}\label{FR2}
\inf_M|R|<\frac{1}{r^2}.
\end{eqnarray}
\end{thm}

\vskip5pt

\noindent{\bf Remark.} For each $a>r$, let $M_a$ be the graph of
$f:\overline B_r\subset\mathbb R^n\to\mathbb R$ given by
$$
f(x_1,...,x_n)=\Big({a^2-\sum_{i=1}^nx_i^2}\Big)^{\frac{1}{2}}.
$$
The mean curvature and the scalar curvature of $M_a$ are,
respectively, $1/a$ and $1/a^2$. Since $a$ can be made arbitrarily
close to $r$, we see that the estimates (\ref{FH1}) and
(\ref{FR2}) are sharp.

\vskip5pt

The key to obtain strict inequalities in (\ref{FH1}) and
(\ref{FR2}) is a general tangency principle by Silva and the
author (\cite{FS1}, Theorem 1.1), which establishes relatively
weak sufficient conditions for two hypersurfaces of a Riemannian
manifold to coincide in a neighborhood of a tangency point (see
Theorems \ref{TPH} and \ref{TPR} in Section 2 for particular cases
of this tangency principle).

\vskip5pt

It is natural trying to establish versions for the higher order
mean curvatures $H_k,\;k\geq 2$ (see Section 2 for the
definitions) of the theorem of Chern and Flanders referred to in
the beginning of the introduction. In this regard, Elbert
\cite{El} proved that there is no entire graph in $\mathbb
R^{n+1}$ with second fundamental form of bounded length and
negative 2-mean curvature $H_2$ bounded away from zero (for
hypersurfaces of a Euclidean space, the 2-mean curvature $H_2$ is
nothing but the scalar curvature $R$ of the hypersurface).
Hasanis-Vlachos \cite{HV} improved Elbert's result by proving that
$\inf_M|R|=0$ for all entire graphs in $\mathbb R^{n+1}$ with
second fundamental form of bounded length (see \cite{El} and
\cite{HV} for results regarding the other higher order mean
curvatures). As an immediate consequence of the first part of
Theorem \ref{FR}, we obtain the following improvement of these
results of Elbert and Hasanis-Vlachos.

\begin{cor}\label{E-HV} If an entire graph $M^n\subset\mathbb R^{n+1}$ has
bounded mean curvature, then
$$
\inf_M|R|=0.
$$
In particular, if the scalar curvature $R$ is constant, then
$R=0$.
\end{cor}

\vskip5pt

A classical theorem by Hilbert states that the hyperbolic plane
cannot be isometrically immersed in the 3-dimensional Euclidean
space. In a remarkable work, Efimov \cite{E} extended Hilbert's
theorem by proving that there is no complete immersed surface in
$\mathbb R^3$ with Gaussian curvature less than a negative
constant.

\vskip5pt

Reilly \cite{R} and Yau \cite{Y1} (see also \cite{Y2}, problem 56,
p. 682) proposed the following extension of Efimov's theorem:

\vskip5pt

``There are no complete hypersurfaces in $\mathbb R^{n+1}$ with
negative Ricci curvature bounded away from zero.''

\vskip5pt

In a well known work, Smith-Xavier \cite{SX} showed that the above
question has a positive answer for $n=3$ and provided a partial
answer for $n>3$. This question has also a positive answer in the
class of all entire graphs with negative Ricci curvature in
Euclidean space, as Chern has shown \cite{C} that
$\inf_M|\text{Ric}|=0$ for all entire graphs $M^n\subset\mathbb
R^{n+1},\;n\geq 3$, with negative Ricci curvature. The corollary
of the following theorem improves this result of Chern.

\vskip10pt

\begin{thm}\label{FRic} Let $M^n\subset\mathbb R^{n+1},\;n\geq 3,$ be a graph over a closed ball $\overline B_r\subset\mathbb R^n$.
If the Ricci curvature of $M$ is negative, then
\begin{eqnarray}\label{FA1}
\inf_M|A|<\frac{3(n-2)}{r},
\end{eqnarray}
where $|A|$ is the length of the second fundamental form of $M^n$
in $\mathbb R^{n+1}$.
\end{thm}

\vskip5pt

\begin{cor}\label{SX} If the Ricci curvature of an entire graph $M^n\subset\mathbb R^{n+1},\,n\geq 3$, is negative,
then {\em inf}$_M|A|=0$.
\end{cor}

Okayasu \cite{O} constructed an example of an $O(2)\times
O(2)$-invariant complete hypersurface of constant negative scalar
curvature in $\mathbb R^4$. Since the length $|A|$ of the second
fundamental form in Okayasu's example is unbounded, one can then
formulate the following Efimov type question: is there a complete
hypersurface in $\mathbb R^{n+1}$ with bounded mean curvature and
negative scalar curvature bounded away from zero? Corollary
\ref{E-HV} shows that if such a hypersurface does exist then
certainly it is not an entire graph. On the other hand, we do not
know whether Corollary \ref{E-HV} holds without the assumption
that the mean curvature is bounded.

\vskip5pt

In dimension 2, Milnor \cite{KO} conjectured (see also \cite{Y2},
problem 62, p. 684) the following improvement of Efimov's result:
If $M^2\subset\mathbb R^3$ is a complete non-flat umbilic free
surface whose Gaussian curvature does not change sign, then $\inf
|A|=0$. Smyth-Xavier \cite{SX} proposed the following analogue in
higher dimensions: If $M^n\subset\mathbb R^{n+1}$ is a complete
immersed hypersurface with negative Ricci curvature, then $\inf_M
|A|=0$. Corollary \ref{SX} shows that this question has a positive
answer for entire graphs in Euclidean spaces.

\vskip10pt

In the following theorem we obtain an estimate for $\inf_M|A|$
under another geometric condition.

\begin{thm}\label{FA} Let $M^n\subset\mathbb R^{n+1}$ be a graph over a closed ball $\overline B_r\subset\mathbb R^n$.
If the mean curvature of $M$ does not change sign, then
\begin{eqnarray}\label{FA2}
\inf_M|A|<\frac{n}{r}.
\end{eqnarray}
\end{thm}

As immediate consequences of Theorem \ref{FA}, we obtain the
following results by Silva and the author \cite{FS2}:

\begin{cor}\label{PH1} If the mean curvature of an entire graph $M^n\subset\mathbb R^{n+1}$ does not change
sign, then {\em inf}$_M|A|=0$.
\end{cor}

Corollary \ref{PH1} was obtained by Hasanis-Vlachos \cite{HV}
under the additional assumption that the length $|A|$ of the
second fundamental form $A$ of $M$ is bounded.

\begin{cor}\label{PH2} Let $M^n\subset\mathbb R^{n+1}$ be an entire graph. If $|A|$ is constant
and $H$ does not change sign, then $M$ is a hyperplane.
\end{cor}

\noindent{\bf Remark.} Corollary \ref{PH1} does not hold for
hypersurfaces which are not graphs. In fact, any circular cylinder
satisfies $\inf |A|>0$.

\vskip10pt

We stress that our methods in this paper are substantially
different from the ones employed by Heinz \cite{H}, which were
based on an ingenious use of the divergence theorem, applied to
the classical formulas for the mean and Gaussian curvature of a
graph in two variables. By contrast, our proofs constitute another
application of our work on the tangency principle \cite{FS1}. They
also use a classical result of G\aa rding \cite{Ga} on hyperbolic
polynomials.

\vskip10pt

\noindent{\bf Acknowledgements.} This paper was partly written
during an extended visit of the author to the University of Notre
Dame. He would like to record his gratitude to the mathematics
department for the invitation and hospitality, and to professor
Frederico Xavier for the suggestions and helpful discussions.

\section{Preliminaries}

Given an oriented immersed hypersurface $M^n$ of the
$(n+1)$-dimensional Euclidean space $\mathbb R^{n+1}$, denote by
$A$ the shape operator associated to the second fundamental form
of the immersion and by $k_1(p),...,k_n(p)$ the principal
curvatures of $M$ at a point $p$, labelled by the condition
$k_1(p)\leq\dots\leq k_n(p)$. The squared length $|A|^2(p)$ of the
second fundamental form at a point $p$ is defined as the trace of
$A^2(p)$. It is easy to see that
\begin{eqnarray}\label{squarelength}
|A|^2(p)=\sum_{i=1}^nk_i^2(p).
\end{eqnarray}
Denote by $R$ the scalar curvature of $M$ and by $H$ the mean
curvature of the immersion. If $e_1,...,e_n$ diagonalizes $A(p)$
with corresponding eigenvalues $k_1,...,k_n$, it follows from the
Gauss equation \cite{Dj} that the Ricci curvature of $M$ at $p$ in
the direction $e_i$ is given by
\begin{eqnarray}\label{RC2}
(n-1)\text{Ric}_p (e_i)=\sum_{j=1 , j\neq i}^n k_i k_j=k_i
(nH-k_i).
\end{eqnarray}
Taking the sum on $i$, we obtain
\begin{eqnarray}\label{HRA}
n^2H^2=|A|^2+n(n-1)R.
\end{eqnarray}

\noindent{For} $1\leq k\leq n$, the $k$-mean curvature $H_k(x)$ of
$M$ at a point $x$ is defined by
\begin{eqnarray}\label{Hr}
H_k(x)=\frac{1}{\binom nk}\sigma_k(k_1(x),...,k_n(x)),
\end{eqnarray}
where $\sigma_k:\mathbb R^n\to\mathbb R$ is given by
\begin{eqnarray}\label{SF}
\sigma_k(x_1,...,x_n)=\sum_{i_1<\dots <i_k}x_{i_1}\dots x_{i_k}
\end{eqnarray}
and is called the $k$-elementary symmetric function. Notice that
$H_1$ is the mean curvature $H$ of the hypersurface and $H_2$ is,
by the Gauss equation \cite{Dj}, simply the scalar curvature $R$
of $M$ (more generally, for hypersurfaces of an ambient space with
constant sectional curvature $c$, we have $R=H_2+c$).

\vskip10pt

For $1\leq k\leq n$, denote by $\Gamma_k$ the connected component
of the set $\{\sigma_k>0\}$ that contains the vector $(1,...,1)$.
It follows immediately from the definitions that $\Gamma_k$
contains the positive cone $\mathcal O^n=\{(x_1,...,x_n)\in\mathbb
R^n:x_i>0,\,\forall i,\}$, for all $1\leq k\leq n$. It was proved
by G\aa rding \cite{Ga} that $\Gamma_k$ is an open convex cone,
$1\leq k\leq n$, and that
\begin{eqnarray}\label{Ga}
\Gamma_1\supset\Gamma_2\supset\dots\supset\Gamma_n.
\end{eqnarray}

Given a hypersurface $M^n\subset\mathbb R^{n+1}$, a point $p\in M$
and a vector $\eta_o\perp T_pM,\;|\eta_o|=1$, we can parametrize a
neighborhood of $p$ in $M$ by
\begin{eqnarray}\label{PG}
\varphi(x)=x+\mu(x)\eta_o,
\end{eqnarray}
for some smooth real valued function $\varphi:V\to\mathbb R$
defined in a neighborhood of 0 in $T_pM$.

\vskip10pt

Let $M_1^n$ and $M_2^n$ be hypersurfaces of $\mathbb R^{n+1}$
tangent at a point $p$ and $\eta_o$ a unitary vector normal to
$T_pM_1=T_pM_2$. Parametrize $M_1$ and $M_2$ as in (\ref{PG}),
obtaining corresponding functions $\mu_1$ and $\mu_2$. As in
\cite{FS1}, we say that $M_1$ remains above $M_2$ in a
neighborhood of $p$ with respect to $\eta_o$ if
$\mu_1(x)\geq\mu_2(x)$ for all $x$ in a neighborhood of zero.

\vskip10pt

In our proofs we will make use of the following theorems, which
are particular cases of a general tangency principle by Silva and
the author (\cite{FS1}, Theorem 1.1).

\begin{thm}\label{TPH}
{\bf (Tangency Principle for Mean Curvature)} Let $M_1^n$ and
$M_2^n$ be hypersurfaces of $\mathbb R^{n+1}$ tangent at a point
$p$ and suppose that $M_1$ remains above $M_2$ in a neighborhood
of $p$ with respect to a unit vector $\eta_o\perp T_pM_1$. If the
mean curvature of $M_2$ at $(x,\varphi_2(x))$ is greater than or
equal to the mean curvature of $M_1$ at $(x,\varphi_1(x))$, for
all $x$ sufficiently small, then $M_1$ and $M_2$ coincide in a
neighborhood of $p$.
\end{thm}

\begin{thm}\label{TPR}
{\bf (Tangency Principle for Scalar Curvature)} Let $M_1^n$ and
$M_2^n$ be hypersurfaces of $\mathbb R^{n+1}$ tangent at a point
$p$ and suppose that $M_1$ remains above $M_2$ in a neighborhood
of $p$ with respect to a unit vector $\eta_o\perp T_pM_1$. If the
scalar curvature of $M_2$ at $(x,\varphi_2(x))$ is greater than or
equal to the scalar curvature of $M_1$ at $(x,\varphi_1(x))$, for
all $x$ sufficiently small, and all principal curvatures
$k_1(p),...,k_n(p)$ of $M_2$ at $p$ are positive (or more
generally, if $\big(k_1(p),...,k_n(p)\big)\in\Gamma_2$), than
$M_1$ and $M_2$ coincide in a neighborhood of $p$.
\end{thm}

\section{Proofs of the Theorems}

\noindent{\bf Proof of Theorem \ref{FH}.} We can suppose
$c:=\text{inf}_M|H|>0$. Otherwise, inequality (\ref{FH1}) is
trivial. Choose the orientation for $M$ so that $H\geq c>0$, and
take a sphere $S$ of radius $r$, disjoint of $M^n$ and contained
in the component of ($\overline B_r\times \mathbb R)\backslash M$
that contains the normals. Move $S$ until it touchs $M$ for the
first time, say at $p$, and denote by $N$ the unit normal vector
field in $M$. By our assumption that $M$ is a graph over
$\overline B_r$ we have that $p$ belongs to the interior of $M$.
If $p_o$ is the center of $S$, we have that $p$ is a point where
the function $f:M\to\mathbb R$, $f(x)=\frac{1}{2}\parallel
x-p_o\parallel^2$, attains its minimum. If $e_1,...,e_n$ is an
orthonormal basis of $T_pM$ such that
$A(e_i)=k_i(p)e_i,\,i=1,...,n,$ we thus have
\begin{eqnarray}\label{grad}
0=\text{grad}f(p)=(p-p_o)^T
\end{eqnarray}
and
\begin{eqnarray}\label{hess}
0\leq \text{Hess}f(p)(e_i,e_i)=1+\langle
\sigma(e_i,e_i),p-p_o\rangle=1+\langle p-p_o,N(p)\rangle k_i(p).
\end{eqnarray}
Equality (\ref{grad}) implies
$$
N(p)=\frac{p_o-p}{\parallel p_o-p\parallel}=\frac{p_o-p}{r}
$$
and, by substitution of this on (\ref{hess}), we conclude that
$k_i(p)\leq\frac{1}{r},\,i=1,...,n$. Thus
\begin{eqnarray}
nc\leq nH(p)=k_1(p)+\dots +k_n(p)\leq \frac{n}{r},
\end{eqnarray}
from which we obtain
\begin{eqnarray}\label{strict1}
\text{inf}_M|H|=c\leq 1/r.
\end{eqnarray}
If equality occurs in (\ref{strict1}), we have $H\geq 1/r$ along
$M$ and, by Theorem \ref{TPH}, $M$ and $S$ coincide in a
neighborhood of $p$. By a connectedness argument, we conclude that
$M$ is a closed hemisphere of $S$. In particular, the tangent
planes to $M$ along $\partial M$ are vertical, contradicting the
assumption that $M$ is a graph over $\overline B_r$. This
contradiction implies that the inequality in (\ref{strict1}) is
strict.\qed

\vskip10pt

\noindent{\bf Proof of Theorem \ref{FR}.} We will first prove
(\ref{FR1}). If $R$ changes sign, there is, by continuity, a point
where the scalar curvature vanishes and (\ref{FR1}) follows
trivially. If $R>0$ along $M$, we have from (\ref{HRA})
\begin{eqnarray}
n(n-1)|R|=n(n-1)R=n^2H^2-|A|^2\leq n^2H^2,
\end{eqnarray}
which implies
\begin{eqnarray}
|R|\leq \frac{nH^2}{n-1}\leq \frac{n}{n-1}|H|\sup|H|.
\end{eqnarray}
Using Theorem \ref{FH}, we obtain
\begin{eqnarray}
\inf_M|R|\leq
\frac{n}{n-1}\sup|H|\inf|H|\leq\frac{n}{r(n-1)}\sup|H|,
\end{eqnarray}
from which we easily obtain (\ref{FR1}).

\vskip5pt

Suppose now $R<0$ everywhere and orient $M$ by a unit normal
vector field $N$. As in the proof of Theorem \ref{FH}, take a
sphere $S$ of radius $r$, disjoint of $M$ and contained in the
component of ($\overline B_r\times \mathbb R)\backslash M$ that
contains the normals, and move $S$ until it touchs $M$ for the
first time, say at $p$. Since $R<0$ along $M$, we have principal
curvatures of both signs at each point of $M$. Let $l$ be number
of negative principal curvatures of $M$ at $p$ so that
\begin{eqnarray}
k_1(p)\leq\dots\leq k_l(p)<0\leq k_{l+1}(p)\leq\dots\leq k_n(p).
\end{eqnarray}
By the Gauss equation, we have
\begin{eqnarray}
0>\frac{n(n-1)}{2}R(p)\nonumber&=&\sum_{1\leq i<j\leq
n}k_ik_j\\&=&\sum_{1\leq i<j\leq l}k_ik_j+\sum_{l+1\leq i<j\leq
n}k_ik_j+\sum_{i=1,...,l;
j=l+1,...,n}k_ik_j\nonumber\\&\geq&\sum_{i=1,...,l;
j=l+1,...,n}k_ik_j,
\end{eqnarray}
and so
\begin{eqnarray}\label{estR}
0>\frac{n(n-1)}{2}R(p)\geq
(k_1+\dots+k_l)(k_{l+1}+\dots+k_n)=\Big(nH-\sum_{i=l+1}^nk_i\Big)\sum_{i=l+1}^nk_i.
\end{eqnarray}
Since $k_i(p)\leq 1/r,\,i=1,\dots,n$ (see the proof of Theorem
\ref{FH}), we arrive at
\begin{eqnarray}
\frac{n(n-1)}{2}\inf|R|&\leq&\nonumber
\frac{n(n-1)}{2}|R(p)|\\\nonumber&\leq&\Big(n\sup|H|+\sum_{i=l+1}^nk_i\Big)\sum_{i=l+1}^nk_i
\\\nonumber&\leq&\Big(n\sup|H|+\frac{n-l}{r}\Big)\frac{n-l}{r}\\&\leq&\Big(n\sup|H|+\frac{n-1}{r}\Big)\frac{n-1}{r},
\end{eqnarray}
from which we easily obtain (\ref{FR1}).

\vskip5pt

We will now proceed to prove the second part of the theorem. Let
$q$ be a point where the second fundamental form is semi-definite
and choose the orientation $N$ so that all principal curvatures at
$q$ are nonnegative. We can suppose $\text{inf}_M|R|>0$, otherwise
there is nothing to prove. Since $k_i(q)\geq 0,\,i=1,...,n$, we
have $R>0$ along $M$. From $\mathcal O^n\subset\Gamma_2$ (see
Section 2) we infer that the principal curvature vector
$\overrightarrow{k}(q)=\big(k_1(q),\dots,k_n(q)\big)$ of $M$ at
$q$ belongs to $\overline{\Gamma_2}$. Since $R(q)>0$, we have in
fact $\overrightarrow{k}(q)\in\Gamma_2$. It follows from the
connectedness of both $M$ and $\Gamma_2$ that
$\overrightarrow{k}(x)\in\Gamma_2$, for all $x\in M$. In
particular, $\overrightarrow{k}(p)\in\Gamma_2$, where $p$ is as in
the first part of the proof . If we had
\begin{eqnarray}\label{infR}
\text{inf}_M|R|\geq 1/r^2,
\end{eqnarray}
we would conclude, by Theorem \ref{TPR}, that $M$ and $S$ coincide
in a neighborhood of $p$. Reasoning as in the proof of Theorem
\ref{FH}, we would conclude that $M$ is a closed hemisphere of
$S$, contradicting the assumption that $M$ is a graph over
$\overline B_r$. This contradiction implies that (\ref{infR}) does
not hold and concludes the proof of the theorem.\qed

\vskip10pt

\noindent{\bf Proof of Theorem \ref{FRic}} Since the Ricci
curvature of $M$ is negative, we have, by (\ref{RC2}), that all
principal curvatures are nonzero and that there are principal
curvatures of both sign at each point of $M$. Let $l$ be the
number of negative principal curvatures, so that $k_1\leq\dots\leq
k_l<0< k_{l+1}\leq\dots\leq k_n$. Since $n\geq 3$, we can choose
the orientation so that $n-1\geq l\geq 2$. Let $S$ and $p$ be as
in the proof of Theorem \ref{FH}, and choose an orthonormal basis
$\{e_1,...,e_n\}$ of $T_pM$ satisfying
$A(e_i)=k_ie_i,\,i=1,...,n$. Since the Ricci curvature is
negative, one has, by (\ref{RC2}),
$$
k_i(k_1+\dots +\widehat{k_i}+\dots +k_l+k_{l+1}+\dots
+k_n)<0,\;\;\;i=1,\dots,l,
$$
and, since $k_i<0$,
$$
k_{l+1}+\dots +k_n>-k_1-\dots -\widehat{k_i}-\dots
-k_l=|k_1|+\dots +\widehat{|k_i|}+\dots +|k_l|,
$$
where the circumflex over $k_i$ means that this term is omitted on
the sum. Taking the sum with $i=1,\dots ,l$, we obtain
\begin{eqnarray}
l(k_{l+1}+\dots +k_n)>(l-1)\sum_{m=1}^l|k_m|.
\end{eqnarray}
Since $k_i(p)\leq 1/r$, $i=1,\dots,n$ (see the proof of Theorem
\ref{FH}), we arrive at
\begin{eqnarray}
\sum_{m=1}^l|k_m|<\frac{l(n-l)}{r(l-1)}.
\end{eqnarray}
Thus
\begin{eqnarray}
\sum_{m=1}^n|k_m|=\sum_{m=1}^l|k_m|+\sum_{m=l+1}^n|k_m|<\frac{l(n-l)}{r(l-1)}
+\frac{n-l}{r}=\frac{(n-l)(2l-1)}{r(l-1)}.
\end{eqnarray}
Noticing that the right hand side of the above equation is
strictly decreasing in $l$, we have
$$
\sum_{m=1}^n|k_m|<\frac{3(n-2)}{r}.
$$
Hence
$$
|A|^2(p)=\sum_{m=1}^n|k_m|^2<\Big(\sum_{m=1}^n|k_m|\Big)^2<\Big(\frac{3(n-2)}{r}\Big)^2,
$$
from which we obtain
$$
\text{inf}|A|\leq |A|(p)<\frac{3(n-2)}{r}.\qed
$$

\vskip10pt

\noindent{\bf Proof of Theorem \ref{FA}.} Choose the orientation
for $M$ so that $H\geq 0$ and let $S$ and $p$ be as in the proof
of Theorem \ref{FH}. We have two cases to consider:

\smallskip

\noindent{\it First case:} All principal curvatures of $M$ at $p$
are nonnegative. Since $k_i(p)\leq 1/r,\,i=1,\dots,n$ (see the
proof of Theorem \ref{FH}), we have
\begin{eqnarray}
|A|^2(p)=\sum_{i=1}^nk_i^2(p)\leq\frac{n}{r^2},
\end{eqnarray}
and so
\begin{eqnarray}
\text{inf}_M|A|\leq |A|(p)\leq\frac{\sqrt n}{r}<\frac{n}{r}.
\end{eqnarray}

\smallskip

\noindent{\it Second case:} There are negative principal
curvatures of $M$ at $p$. Let $l$ the number of negative principal
curvatures so that
\begin{eqnarray}
k_1(p)\leq\dots\leq k_l(p)<0\leq k_{l+1}(p)\leq\dots\leq k_n(p).
\end{eqnarray}
Notice that $l\leq n-1$ since $H\geq 0$. From $k_i(p)\leq
1/r,\,i=1,\dots,n$, and $H\geq 0$, we obtain
\begin{eqnarray}
\frac{n-l}{r}\geq k_{l+1}(p)+\dots+k_n(p)\geq -k_1(p)-\dots
-k_l(p)=|k_1|(p)+\dots +|k_l|(p).
\end{eqnarray}
Hence,
\begin{eqnarray}
|A|^2(p)&=&\sum_{i=1}^lk_i^2+\sum_{i=l+1}^nk_i^2\leq\Big(\sum_{i=1}^l|k_i|\Big)^2+\sum_{i=l+1}^nk_i^2\nonumber
\\&\leq&
\frac{(n-l)^2}{r^2}+\frac{n-l}{r^2}=\frac{(n-l)(n-l+1)}{r^2}\nonumber\\&\leq&\frac{n(n-1)}{r^2}<\frac{n^2}{r^2},
\end{eqnarray}
from which we obtain (\ref{FA2}). \qed

\noindent{Francisco Fontenele\\Departamento de Geometria\\
Universidade Federal Fluminense\\Niter\'oi, RJ, Brazil\\e-mail:
fontenele@mat.uff.br}


\begin{thebibliography}{s2}

\bibitem{BBM} J. L. Barbosa, G. P. Bessa and J. F. Montenegro, {\em On Bernstein-Heinz-Chern-Flanders inequalities},
Math. Proc. Cambridge Philos. Soc., {\bf{144}} (2008) 457-464.

\bibitem{C} S. S. Chern, {\em On the curvatures of a piece of hypersurface in Euclidean space},
Abh. Math. Sem. Univ. Hamburg, {\bf{29}} (1965) 77-91.

\bibitem{Dj} M. Dajczer, {\em Submanifolds and isometric immersions},
Mathematics Lectures Series, 13, Houston: Publish or Perish, 1990.

\bibitem{E} N. V. Efimov, {\em Hyperbolic problems in theory of
surfaces}, Proc. Int. Congress Math. Moscow (1966); Amer. Math.
Soc. Transl., {\bf{70}} (1968) 26-38.

\bibitem{El} M. F. Elbert, {\em On complete graphs with negative r-mean curvature},
Proc. Amer. Math. Soc., {\bf{128}} (2000) 1443-1450.

\bibitem{Fl} H. Flanders, {\em Remark on mean curvature},
J. London Math. Soc., {\bf{41}} (1966) 364-366.

\bibitem{FS1} F. Fontenele and S. L. Silva, {\em A tangency principle and applications},
Illinois J. Math., {\bf{45}} (2001) 213-228.

\bibitem{FS2} F. Fontenele and S. L. Silva, {\em Sharp estimates for the size of balls in the complement of a hypersurface},
Geom. Dedicata, {\bf{115}} (2005) 163-179.

\bibitem{Ga} L. G\aa rding, {\em An inequality for hyperbolic polynomials},
J. Math. Mech., {\bf{8}} (1959) 957-965.

\bibitem{HV} T. Hasanis and T. Vlachos, {\em Curvature properties of hypersurfaces},
Arch. Math.(Basel), {\bf{82}} (2004) 570-576.

\bibitem{H} E. Heinz, {\em \"{U}ber Fl\"{a}chen mit eineindeutiger Projektion auf eine Ebene,
deren Kr\"{u}mmungen durch Ungleichungen eingeschr\"{a}nkt sind.},
Math. Ann., {\bf{129}} (1955) 451-454.

\bibitem{KO} T. Klotz and R. Osserman, {\em Complete surfaces in $E^3$ with constant mean curvature},
Comment. Math. Helv., {\bf{41}} (1966-67) 313-318.

\bibitem{O} T. Okayasu, {\em $O(2)\times O(2)$-Invariant hypersurfaces with constant negative scalar curvature in $E^4$},
Proc. Amer. Math. Soc., {\bf{107}} (1989) 1045-1050.

\bibitem{R} R. Reilly, {\em On the Hessian of a function and the curvatures of its graph},
Michigan Math. J., {\bf{20}} (1973) 373-383.

\bibitem{Sa} I. Salavessa, {\em Graphs with parallel mean curvature},
Proc. Amer. Math. Soc., {\bf{107}} (1989) 449-458.

\bibitem{SX} B. Smith and F. Xavier, {\em Efimov's theorem in dimension greater than two},
Invent. Math., {\bf{90}} (1987) 443-450.

\bibitem{Y1} S. T. Yau, {\em Submanifolds with constant mean curvature II},
Amer. J. Math., {\bf{97}} (1975) 76-100.

\bibitem{Y2} S. T. Yau, {\em Seminar on differential geometry},
Annals of Mathematics Studies 102. Princeton, N.J.: Princeton
Univ. Press, 1982.



\end{thebibliography}
\end{document}